\documentclass[final,authoryear]{elsarticle}

\journal{MTNS 2024}

\def\appendixname{}
\usepackage{graphicx}
\usepackage{subcaption}
\usepackage{comment}
\usepackage{tikz}
\usetikzlibrary{calc}
\usepackage{color}
\usepackage{url}
\usepackage{amsthm,amsmath,yhmath,amssymb,amsfonts,mathtools,mathrsfs}
\usepackage{natbib}

\excludecomment{verMTNS}
\includecomment{verArxiv}

        {\hspace*{\fill}$\Box$\par}

\theoremstyle{plain}
\newtheorem{thm}{Theorem}
\newtheorem{lem}[thm]{Lemma}
\newtheorem{prop}[thm]{Proposition}

\theoremstyle{definition}
\newtheorem{defn}[thm]{Definition}

\usepackage{defs}

\begin{document}

\begin{frontmatter}

\title{Convergence Analysis of Ensemble Filters for Linear Stochastic Systems with Poisson-Sampled Observations\tnoteref{anr}}

\tnotetext[anr]{The work of A.~Tanwani is partly supported by the project {\sc CyPhAI}, financed by ANR-JST CREST program with grant number ANR-20-JSTM-0001. 
The work of O. Yufereva  was performed as a part of research conducted in the Ural Mathematical Center with the financial support
of the Ministry of Science and Higher Education of the Russian
Federation (Agreement number 075-02-2024-1377).
}

\author[First]{Aneel Tanwani}
\author[First,Second]{Olga Yufereva}

\address[First]{LAAS -- CNRS, University of Toulouse, CNRS, Toulouse France.\\
Email: {\tt aneel.tanwani@cnrs.fr}}
\address[Second]{Krasovskii Institute of Mathematics and Mechanics, Ural Branch of RAS, Yekaterinburg, Russia. 
Email: {\tt oiufereva@laas.fr}}


\begin{abstract}
For continuous-time linear stochastic dynamical systems driven by Wiener processes, we consider the problem of designing ensemble  filters when the observation process is randomly time-sampled. We propose a continuous-discrete McKean--Vlasov type diffusion process with additive Gaussian noise in observation model, which is used to describe the evolution of the individual particles in the ensemble. These particles are coupled through the empirical covariance and require less computations for implementation than the optimal ones based on solving Riccati differential equations. Using appropriate analysis tools, we show that the empirical mean and the sample covariance of the ensemble filter converges to the mean and covariance of the optimal filter if the mean sampling rate of the observation process satisfies certain bounds and as the number of particles tends to infinity.
\end{abstract}


\end{frontmatter}

\section{Introduction}

A fundamental problem of interest in modern engineering systems is to develop techniques for extracting information about the unknown elements of the system from the data measurements. Looking at such problems from the viewpoint of dynamical systems, we can naturally study them using filtering algorithms. 
The early developments in the area of filtering can be traced back to \citep{ref:KalBuc-61}, where the authors provided an elegant recursive closed form solution for optimal state estimation in linear stochastic systems. Since this pioneering work, the field has evolved enormously to develop algorithms for nonlinear dynamical systems, address different information structures, and provide computationally efficient methods for large-scale systems. An excellent compilation of different developments on the topic of filtering appears in \citep{ref:CriRoz-11}.

The basic problem in optimal filtering is to compute the distribution of the unknown state process conditioned upon the observations; see \citep{ref:AndMoo-69} for an overview.
Among the existing techniques for filtering, the use of Monte Carlo integration methods for approximating the optimal distribution have gained significant interest in the literature \citep{doucet2001sequential,ristic2003beyond}. In the same spirit, \cite{ref:Eve-94} introduced the technique of ensemble Kalman filters, or feedback particle filters, to develop filtering methods for large scale applications related to geophysical sciences, so that the error covariance is computed from a {\em collection} of state estimators rather than from a single (Riccati) differential equation. Since then, the use of ensemble Kalman filters has had a notable impact where estimation with noisy data is required in large-scale models. This approach is based on simulating the evolution of individual particles in the ensemble through differential equations that are coupled through the associated empirical mean and the empirical error covariance. Several review articles \citep{chen2003bayesian, ref:CavReiStu-22} and the books \citep{ref:Eve-09} provide an overview of developments in that area. In most of these works, we do not find much details about the theoretical analysis of the proposed filtering techniques, and this area of mathematical analysis of the ensemble Kalman filters has gathered attention only very recently.

 Over the past decade, the ensemble filtering methods have been viewed from the lens of mean-field models described by stochastic differential equations, and the ensemble particles are simply the approximations of these mean field models. Recent review articles which elaborate on this viewpoint are \citep{bishop2020mathematical,ref:TagMeh-23}.
In fact, the limiting behavior of these particles is described by a McKean--Vlasov type diffusion process, which is also referred to as the mean-field process. 
In the literature, this mean-field process is chosen in different ways, e.g., by adding noise in the prediction term and the correction term of the Kalman--Bucy process \citep{ref:DelTug-18}, or as a non-diffusion equation that is optimal in the measure transportation sense \citep{ref:TagMeh-20}.

Another important research direction is to study the filtering problem with constraints on the information available for computing the optimal distribution. In particular, for implementation of filters subject to observations transmitted through some communication protocols, it is natural to stipulate that the observations arrive at some random time instants \citep{ref:SinSch-04}. It is of interest to compute the conditional distribution of the state process conditioned upon this discrete observation process~\citep{ref:Jaz-07}, and results in continuous-discrete filters. 
For certain technical reasons, and to better study the effect of mean sampling rate, we stipulate in our previous works that the sampling process is a Poisson counter. In particular, for a system class very close to the one studied in this paper, \citep{ref:TanYuf-20} proposes a continuous-discrete filter and analyze the boundedness of error covariance as a function of the mean sampling rate. A similar performance is studied for a collection of interacting filters over a network in \citep{ref:Tan-22}.

Our primary objective in this article is to develop ensemble filters for continuous-time stochastic processes subject to randomly time-sampled observations. In the literature, we find some variants of {\em continuous-discrete feedback particle filters} in different settings. The paper \citep{ref:BerRei-10} provides one (and possibly the first) such example, where the authors use mollifiers in the particle equations to smoothen the dynamics, but no statements about the limiting process are provided. The paper \citep{ref:YanMeh-13} develops particle filters for nonlinear systems using the time-discretization procedure as a part of the derivation and studies convergence as the length of the sampling interval converges to zero. 
In our recent work \citep{ref:YufTan-23}, we developed the transport-inspired particle filters with randomly sampled observations

For our purposes, the state process is modeled by linear continuous-time Ornstein-Uhlenbeck process and the sampling process for the observations is a Poisson counter. We develop ensemble filters which update their estimate whenever the Poisson counter increments due to the arrival of a new measurement from the observation process. 
In contrast to \citep{ref:YufTan-23}, our objective here is to study a different class of ensemble filters which is basically a generalization of the so-called {\em vanilla filters}: It involves noise terms in the description of the mean-field process. The presence of noise is motivated by the fact that it leads to positive definite (and hence invertible) error covariance matrix that is important for good performance of the ensemble particles. However, the presence of noise terms makes the analysis more challenging, which are tackled using  appropriate tools. In essence, we show that for appropriate sampling rate, if the number of the particles in the ensemble is large enough then the empirical mean and the empirical covariance converge to the optimal mean and the covariance.

\section{Problem Setup}\label{sec:prob}

Let us begin with the description of the system class and the formulation of the basic filtering problem.

\subsection{System Class}
We consider dynamical systems modeled by linear stochastic differential equations of the form
\begin{equation}\label{eq:stateProc}
\dd x_t  = A x_t \, \dd t + G\, \dd \omega_t
\end{equation}
where $\proc{x_t}_{t \ge 0}$ is an $\R^n$-valued diffusion process describing the state. Let $(\Omega,\cF, \PP)$ denote the underlying probability space. It is assumed that, for each $t \ge 0$, $\proc{\omega_t}_{t\ge 0}$ is a zero mean $\R^m$-valued standard Wiener process adapted to the filtration $\cF_t \subset \cF$, with the property that $\EE\expecof{\dd\omega_t \, \dd\omega_t^\top} = I_{m} \dd t$, for each $t \ge 0$. The matrices $A \in \R^{n \times n}$ and $G \in \R^{n\times m}$ are taken as constant with $(A,G)$ controllable, and the process $\proc{\omega_t}_{t\ge 0}$ does not depend on the state. The solutions of the stochastic differential equation \eqref{eq:stateProc} are interpreted in the sense of It\^o stochastic integral.

\subsection{Measurement process}

Our goal is to study the state estimation problem when the output measurements are available only at random times. The motivation to work with randomly time-sampled measurements comes from several applications, such as, communication over networks which allow information packets to be sent at some discrete randomly distributed time instants.
Thus, we consider a monotone nondecreasing sequence \(\proc{\tau_n}_{n\in\Nz}\) taking values in \(\R_{\ge 0}\) which denote the time instants at which the measurements are available for estimation.
We introduce the process $N_t$ defined as
\begin{equation}\label{eq:defNt}
	N_t \Let \sup\set[\big]{n\in\Nz \suchthat \tau_n \le t}\quad \text{for }t \in \R,
\end{equation}
and it is assumed that $\proc{N_t}_{t \ge 0}$ is a Poisson process of intensity $\lambda > 0$ and it is independent of the noise and the state processes. Recall \cite[Theorem 2.3.2]{ref:SuhKel-08} that the Poisson process\index{Poisson process} of intensity \(\lambda > 0\) is a continuous-time random process \(\proc[\big]{N_t}_{t\ge 0}\) taking values in \(\N := \Nz \cup \{0\}\), with \(N_0 = 0\), for every \(n\in\N\) and \(0 \teL t_0 < t_1 < \cdots < t_n < +\infty\), the increments \(\{N_{t_k} - N_{t_{k-1}}\}_{k=1}^n\) are independent, and \(N_{t_k} - N_{t_{k-1}}\) is distributed as a Poisson-\(\lambda(t_k - t_{k-1})\) random variable for each \(k\). The Poisson process is among the most well-studied processes, and standard results (see, e.g., \cite[\S2.3]{ref:SuhKel-08}) show that it is memoryless and Markovian.

The discretized, and noisy, observation process is thus defined as
\begin{equation}\label{eq:discOutput}
y_{\tau_{N_t}} = C x(\tau_{N_t}) + \nu_{\tau_{N_t}}, \qquad t \ge 0,
\end{equation}
where $C \in \R^{p \times n}$ is a constant matrix, and $\nu_{k}$ is a sequence of i.i.d. Gaussian noise processes and $\nu_0 \sim \cN(0,V)$.
Equation~\eqref{eq:discOutput} is motivated by the fact that a continuous observation process $\dd z = Cx \dd t + d\eta$ for a Wiener process $\proc{\eta_t}_{t\ge 0}$ is formally equivalent to $y_t = Cx_t + \nu_t$, with the identifications $y_t \sim \frac{\dd z_t}{\dd t}$ and $\nu_t \sim \frac{\dd \eta_t}{\dd t}$, so that $\nu_t$ is a  Gaussian process; see \cite[Chapter~4]{ref:Jaz-07} for further details.
Our goal is to construct the estimate $\widehat x_t$, which minimizes the mean square estimation error, using the observations $\cY_t:=\{y_{\tau_k} \, \vert \, k \le N_t \}$.

\subsection{Optimal filter}
The basic problem in filter design is to find an estimate of the state process which minimizes the mean square estimation error, and is described by the expectation of the  state process $\proc{x_t}_{t\ge 0}$ conditioned upon the measurements observed over the interval $[0,t]$, that is, $\cY_t$.  In particular, with the structure imposed on the system dynamics in this section, the conditional expectation is Gaussian and the two moments are simulated through ordinary differential equations with updates at times when a new measurement arrives. 
For an arbitrary strictly increasing real-valued sequence $(\tau_k)_{k\in \N}$, this procedure is also proposed in \cite[Thm.~7.1]{ref:Jaz-07}. If we specify a sequence $(\tau_k)_{k\in \N}$ so that it corresponds to the arrival times of a Poisson process, we simulate the mean of the conditional distribution as:
\begin{subequations}\label{eq:optMean}
\begin{align}
\dot{\wh X_t} &= A \wh X_t \dd t , \quad t \in [\tau_{N_t}, \tau_{1+N_t}[ \label{eq:optMeana}\\
\wh X_t^+ & = \wh X_t + K_t (y_t - C \wh X_t), \quad t = \tau_{N_t} \label{eq:optMeanb}
\end{align}
\end{subequations}
where the injection gain $K_t = P_t C^\top (C P_t C^\top + V_t)^{-1}$, and the error covariance process $\proc{P_t}_{t \ge 0}$ is described as
\begin{subequations}\label{eq:optVar}
\begin{align}
\dot {P_t} &= (A P_t + P_t A^\top + G G^\top)  , \quad t \in [\tau_{N_t}, \tau_{1+N_t}[ \label{eq:optVara} \\
\wh P_t^+ & = P_t - P_t C^\top (C P_t C^\top + V_t)^{-1} C P_t , \quad t = \tau_{N_t}. \label{eq:optVarb}
\end{align}
\end{subequations}
To write things more compactly later on, we adopt the formalism of writing the continuous-discrete equations \eqref{eq:optMeana} and \eqref{eq:optMeanb} together in a single differential equation, when the jumps are driven by a Poisson counter $N_t$ :
\begin{equation}\label{eq:optMeanPoisson}
\dd \wh X_t = A \wh X_t \dd t + K_t (y_t - C \wh X_t) \dd N_t.
\end{equation}
Similarly, using this formalism, \eqref{eq:optVara} and \eqref{eq:optVarb} can be written in a combined form as,
\begin{equation}\label{eq:optVarPoisson}
\dd P_t = (A P_t + P_t A^\top + G G^\top) \dd t  - P_t C^\top (C P_t C^\top + V_t)^{-1} C P_t \dd N_t.
\end{equation}
In the foregoing discussion, one makes the observation that the optimal conditional distribution is Gaussian for each realization of $(N_t)_{t \ge 0}$ despite the fact that the mean and covariance are discontinuous along each sample path.

From analysis viewpoint, it is important to look at the expectation of the process $\proc{P_t}_{t\ge 0}$ with respect to the sampling times $\proc{\tau_{N_t}}_{t\ge 0}$. In our previous work \citep{ref:TanYuf-20}, we showed that the expectation of piecewise deterministic process $P_t$, denoted by $\cP_t$, is described by the following differential equation:
\begin{equation}
\dot \cP_t = A \cP_t + \cP_t A^\top + GG^\top - \lambda \cP_t C^\top (C \cP_t C^\top + V)^{-1} C\cP_t.
\end{equation}
We can now provide conditions in terms of the bounds on the mean sampling rate $\lambda > 0$ and the structural assumptions on controllability and observability of the pairs $(A,G)$ and $(A,C)$ that guarantee boundedness of $\cP_t$. The boundedness is also important for asymptotic analysis of the first moment of the error process $\proc{x_t - \wh X_t}_{t \ge 0}$.

\section{Diffusion Process and Ensemble Filters}\label{sec:particles}

In the previous section, we saw that the implementation of conventional optimal filter, even for linear systems, is computationally heavy since it requires simulating a Riccati differential equation to compute the injection gains. For a state process evolving in $\R^n$, this optimal filter involves solving $(n^2+n)/2$ differential equations for the covariance process $\proc{P_t}_{t \ge 0}$ and $n$ differential equations for the first moment of the estimate $\proc{\wh X_t}_{t \ge 0}$. This can be quite cumbersome for large values of $n$ especially when we take into consideration the additional operations involved in computing these processes. For this reason, there has been extensive research for other methods to implement optimal filter, or its approximation. One such technique is based on the use of the {\em ensemble}, or feedback particle, filters. The basic idea of the particle filters is to simulate a collection of particles through stochastic differential equation which are coupled to each other through joint statistics of the population.

\subsection{Description of McKean-Vlasov type process}
For the definition of particles in the ensemble, we first describe a diffusion process of McKean-Vlasov type differential equation. This process is, in particular, continuous-discrete because our observation process is discrete. 
\begin{subequations}\label{eq:vanillaMcKean}
\begin{align}
    \dd S_t &:= A S_t \dd t + G \dd \bar \omega_t + Q_t C^\top(CQC^\top + V_t)^{-1} \left[ y_t - (CS_t + \bar \nu_t) \right] \dd N_t \label{eq:vanillaMcKeana} \\
    Q_t &:= \EE \cexpecof{(S_t - \wh S_t) (S_t - \wh S_t)^\top  \given \cY_t } \label{eq:vanillaMcKeanb} \\
        \wh S_t &:= \EE[S_t\mid \cY_t], \label{eq:vanillaMcKeanc}
\end{align}
\end{subequations}
where the processes $\bar \omega_t$ and $\bar \nu_t$ are independent copies of the processes $\omega_t$ and $\nu_t$ respectively, that were introduced in \eqref{eq:stateProc} and \eqref{eq:discOutput}. The motivation of working with this particular diffusion process comes from the fact that the conditional mean $\wh S_t$ and the conditional covariance $Q_t$ satisfy the same differential equations as the ones obtained for their optimal counterparts, that is, it can be checked that
\begin{equation}
\dd \wh S_t = A \wh S_t \dd t + Q_t C^\top(CQ_tC^\top + V_t)^{-1}(y_t - C\wh S_t) \dd N_t
\end{equation}
and
\begin{equation}
\begin{aligned}
\dd Q_t & = ( A Q_t + Q_t A^\top + GG^\top ) \dd t - Q_t C^\top(CQ_tC^\top + V_t)^{-1} C Q_t \dd N_t,
\end{aligned}
\end{equation}
so that, one sees the resemblance with the processes described in \eqref{eq:optMeanPoisson} and \eqref{eq:optVarPoisson}, respectively.

\subsection{Ensemble Filters}
From the mean-field process \eqref{eq:vanillaMcKean}, we describe an ensemble of $M$ particles. For $\ell = 1, \dots, M,$ the particle $S_t^\ell$ is described by the following stochastic differential equation conditioned upon the observation process and the Poisson sampling process:
\begin{equation}\label{eq:defParticles}
\dd S_t^\ell = A S_t^\ell \dd t + G \dd\omega_t^\ell + K_t\upN(y_t - CS_t^\ell - \nu_t^\ell) \, \dd N_t,
\end{equation}
where, for each $\ell$, $\omega_t^\ell$ and $\nu_t^\ell$ represent the independent copies of the processes $\omega_t$ and $\nu_t$ respectively. The term $K\upN_t$ denotes the empirical gain computed at time instants $t = \tau_{N_t}$ as follows:
\begin{equation}\label{eq:defEmpGain}
K\upN_t = Q\upN_t C^\top (C Q\upN_t C^\top + V)^{-1} \\
\end{equation}
in which $Q\upN_t$ denotes the empirical covariance, which is computed using the empirical mean $\wh S\upN_t$ as follows:
\begin{align}
Q\upN_t & = \frac{1}{\pNum} \sum_{\ell = 1}^{\pNum} (S_t^\ell - \wh S\upN_t)(S_t^\ell - \wh S\upN_t)^\top \label{eq:empVarParticles}\\
\wh S\upN_t & = \frac{1}{\pNum} \sum_{\ell = 1}^{\pNum} S_t^\ell. \label{eq:empMeanParticles}
\end{align}
One clearly sees that the computational effort in implementing the ensemble filters \eqref{eq:defParticles} involves solving $\pNum n$ differential equations. For computing the gain, the empirical mean and the sample covariance only needs to be updated at discrete times using relatively simpler arithmetic operations.

\subsection{Convergence result}
The obvious question concerning the proposed filters in \eqref{eq:defParticles} is whether they approximate the optimal filter \eqref{eq:optMeanPoisson}--\eqref{eq:optVarPoisson} reasonably well. In this regard, our main result shows that if the mean sampling rate satisfies certain bounds, then the difference (in expectation) between the empirical covariance of the particles \eqref{eq:empVarParticles} and the optimal covariance \eqref{eq:empMeanParticles} is $O(\frac{1}{M})$. Similarly, we can show that, as the number of particles increase, the expectation of the empirical mean of the particles approaches the expected value of the optimal estimate $\wh X$. 

Our main result formalizing these observations is stated next. To avoid technicalities in the proof, we present the convergence result for the scalar case only, where $A$ is a real number and $G, C$ are nonzero scalars. This allows us to give explicit bounds on the approximation error, and the mean sampling rate required for establishing these bounds.
\begin{thm}\label{thm:main}
Consider system~\eqref{eq:stateProc} with optimal estimator \eqref{eq:optMeanPoisson} and the ensemble of $M$ filters described by \eqref{eq:defParticles}.
Let 
$
\ol \gamma := \max\left\{\frac{V^2}{(\min\{P_0,Q_0^M\} C^2 + V)^2}, \frac{V^2 (\lambda - 2A)^2}{((1-1/M)G^2 C^2 + V (\lambda - 2A))^2} \right\}.
$ 
If the mean sampling rate satisfies the inequality 
\begin{equation}\label{eq:lowerBndLambda}
\lambda > 2A/(1-\ol \gamma^2)
\end{equation}
then it holds that,
\begin{equation}\label{eq:ultBndErrorApproxVar}
\limsup_{t\to\infty} \| \EE\expecof{P_t} - \EE\expecof{Q\upN_t} \| \le \frac{G^2 + \lambda V/C^2 } {M(\lambda(1-\ol\gamma^2) - 2 A)}
\end{equation}
and moreover, 
\[
\| \EE\expecof{ \wh S\upN_t } - \EE\cexpecof{\wh X_t } \| \rightarrow 0, \quad \text{ as } \pNum \rightarrow \infty, \text {and } t \rightarrow \infty. 
\]
\end{thm}
The proof of this result is carried out in Section~\ref{sec:analysis}. The constant $\ol\gamma$ appearing in the theorem statement actually corresponds to the lower bound on the covariance matrices $P_t$ and $Q\upN_t$. Nonzero value of $G$ is important for establishing this lower bound. Because of nonlinear dependence of $\ol\gamma$ on $\lambda$, inequality \eqref{eq:lowerBndLambda} provides the range of values for the mean sampling rate for which the result in Theorem~\ref{thm:main} holds. For the measurement covariance $V$ small enough, this inequality is indeed feasible. Also, if $A \le 0$, \eqref{eq:lowerBndlambda} trivially holds with any $\lambda > 0$.
The ultimate bound on approximation of covariance matrix in \eqref{eq:ultBndErrorApproxVar} increases with the increase in the covariance of process and sensor noises, and decreases as number of particles $M$ gets large.

\subsection{Simulation of an academic example}
To illustrate the simulation of our proposed ensemble filter in \eqref{eq:defParticles} and the result reported in Theorem~\ref{thm:main}, we consider an academic example where the system is described by the equations:
\begin{subequations}\label{eq:linSys}
\begin{align}
\dd x_t &= A x_t \dd t + G \dd \omega_t \label{eq:linSysa}\\
y_{\tau_k} &= C x_{\tau_k}  +  \nu_{\tau_k}, \label{eq:linSysb}
\end{align}
\end{subequations}
with $A = \begin{bsmallmatrix} 0 & 3 &1 \\ 2 & -2 & 1 \\ -2 & 1 & -3 \end{bsmallmatrix}$, $C = \begin{bsmallmatrix}1 & -1 & 2\\ 1 & 0 & 1\end{bsmallmatrix}$, $G = \begin{bsmallmatrix}0.5 & 0.5 & 0.5\end{bsmallmatrix}^\top$, and for each $k \in \N$, $\nu_{\tau_k}$ is normally distributed with mean $(0,0)^\top$ and the constant variance $V = \begin{bsmallmatrix}0.5 & 0.1 \\ 0.1 & 0.5\end{bsmallmatrix}$.

\begin{figure}[!h]
\centering
\begin{subfigure}[t]{0.48\linewidth}\centering
\includegraphics[width=0.99\linewidth]{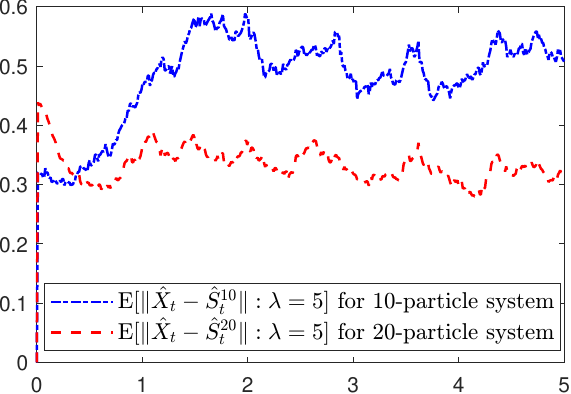}
\caption{The difference of optimal mean and empirical mean is plotted for 10-particles (blue curve) and 20 particles (red curve).}
\label{fig:one_path-a}
\end{subfigure}
\begin{subfigure}[t]{0.48\linewidth}\centering
\includegraphics[width=0.99\linewidth]{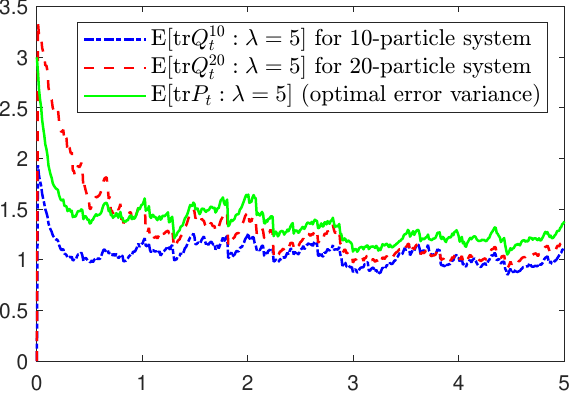}
\caption{The plot shows the optimal covariance (green); empirical covariance for 10-particles (blue) and 20 particles (red).}
\label{fig:smoothed_paths-a}
\end{subfigure}
\caption{Comparison of mean and covariance for the optimal and ensemble filters with mean sampling rate $\lambda = 5$.}
\label{fig:one}
\end{figure}

\begin{figure}[!h]
\centering
\begin{subfigure}[t]{0.48\linewidth}\centering
\includegraphics[width=0.99\linewidth]{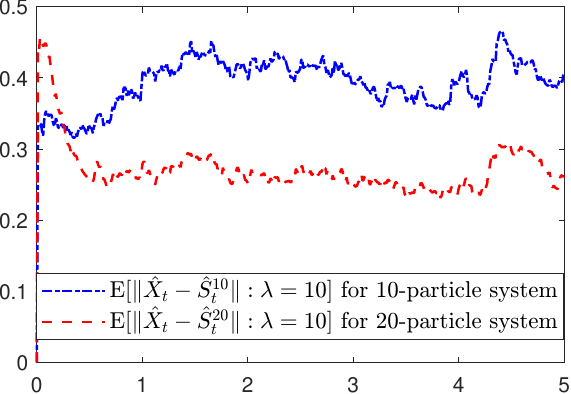}
\caption{The difference of optimal mean and empirical mean is plotted for 10-particles (blue curve) and 20 particles (red curve).}
\label{fig:one_path-b}
\end{subfigure}
\begin{subfigure}[t]{0.48\linewidth}\centering
\includegraphics[width=0.99\linewidth]{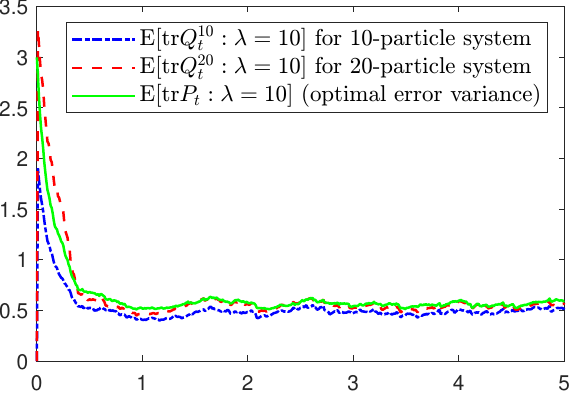}
\caption{The plot shows the optimal covariance (green); empirical covariance for 10-particles (blue) and 20 particles (red).}
\label{fig:smoothed_paths-b}
\end{subfigure}
\caption{Comparison of mean and covariance for the optimal and ensemble filters with mean sampling rate $\lambda = 10$.}
\label{fig:two}
\end{figure}

To measure the effectiveness of the ensemble filter, we compare it with the optimal estimator, and look at the entities that are being compared in Theorem~\ref{thm:main} for different values of mean sampling rate. The simulation results reported in Figure~\ref{fig:one} and Figure~\ref{fig:two} are consistent with the results in Theorem~\ref{thm:main}.


\section{Convergence of Ensemble Filters}\label{sec:analysis}
\begin{verMTNS}
This section essentially carries out the proof of Theorem~\ref{thm:main}. We first derive the differential equations for empirical mean and covariance using Ito's chain rule. We then compute the expectations of the processes using the generator equation, which are then used for studying convergence results stated in Theorem~\ref{thm:main}. For the sake of completeness, these results are stated in the extended version \citep{ref:TanYuf-24}. We emphasize that the developments in this section are carried out for the scalar case only.
\end{verMTNS}

\begin{verArxiv}
This section essentially carries out the proof of Theorem~\ref{thm:main}. We first derive the differential equations for empirical mean and covariance using Ito's chain rule, using Proposition~\ref{prop:chainRuleIto} in Appendix~\ref{sec:tools}). We then compute the expectations of the processes using the generator equation stated in Proposition~\ref{prop:generator}, which are then used for studying convergence results stated in Theorem~\ref{thm:main}. We emphasize that the developments in this section are carried out for the scalar case only.
\end{verArxiv}

\subsection{Evolution of empirical mean and covariance}
We recall that the evolution of the particles is described by the differential equation~\eqref{eq:defParticles} and the corresponding empirical mean is defined in \eqref{eq:empMeanParticles}. We now write down the evolution of the empirical mean of the particles in differential form as follows:
\[
\begin{aligned}
\dd \wh S\upN_t & = A \wh S\upN_t + \frac{1}{\sqrt{\pNum}} G \dd \wt \omega\upN_t + K\upN_t \left( y_t -C\wh S\upN_t - \frac{1}{\sqrt{M}} \wt \nu\upN_t \right) \dd N_t
\end{aligned}
\]
where we used the notation
\[
\wt \omega\upN_t : = \frac{1}{\sqrt{\pNum}} \sum_{\ell = 1}^{\pNum} \omega_t^\ell \, \quad \text { so that } \quad \EE\expecof{\dd\wt \omega\upN_t \dd {\wt \omega^{\pNum \top} _t}} = \dd t
\]
and 
\[
\wt \nu\upN_t := \frac{1}{\sqrt{\pNum}} \sum_{\ell = 1}^{\pNum} \nu_t^\ell \, \quad \text { so that } \quad \wt \nu\upN_t \sim \nu.
\]
Note that the process $\wt \nu\upN_t$ is also i.i.d. Gaussian with covariance equal to $V$. 

Next we write the empirical covariance of the particles, defined in \eqref{eq:empVarParticles}, in differential form. To do so, we let
\[
q_t^\ell := S_t^\ell - \wh S\upN_t.
\]
It then follows that 
\begin{equation}
\dd q_t^\ell = A q_t^\ell \dd t+ G \left( \dd \omega_t^\ell - \frac{1}{\sqrt{M}} \dd \wt \omega\upN_t\right)  - K\upN_t C q^\ell_t \dd N_t -  K\upN_t \left( \nu_t^\ell - \frac{1}{\sqrt{M}} \wt \nu\upN_t \right) \dd N_t.
\end{equation}
\begin{verArxiv}
Consequently, we can use Proposition~\ref{prop:chainRuleIto} to get    
\end{verArxiv}
\begin{verMTNS}
Consequently, we can use Ito's chain rule to get        
\end{verMTNS}
\[
\begin{aligned}
& \dd (q_t^\ell)^2 =  \left (2 A (q_t^\ell)^2 + G^2 - \frac{1}{M} G^2\right) \dd t 
+ 2 G q_t^\ell \left(\dd \omega_t^\ell - \frac{1}{\sqrt{M}}\dd \wt \omega\upN_t \right) \\
& \qquad \quad + \Big[\big(q_t^\ell - K\upN_t C q_t^\ell - K\upN_t \nu_t^\ell +\frac{1}{\sqrt{M}}K\upN_t \wt \nu\upN_t \big)^2 - (q_t^\ell)^2 \Big] \dd N_t.
\end{aligned}
\]
Straightforward simplification yields
\[
\begin{aligned}
\dd (q_t^\ell)^2 & =  \left (2 A (q_t^\ell)^2 + G^2 - \frac{1}{M} G^2\right) \dd t 
+ 2 G q_t^\ell \left(\dd \omega_t^\ell - \frac{1}{\sqrt{M}}\dd \wt \omega\upN_t \right) \\
& \quad+ \Big[ (K^{\pNum}_t)^2 C^2 (q_t^\ell)^2- 2 K\upN_t C (q_t^\ell)^2 + (K^{\pNum}_t)^2 (\nu_t^\ell)^2 \\ 
& \quad + \frac{1}{M}(K^{\pNum}_t)^2 (\wt \nu\upN_t)^2 - \frac{2}{\sqrt{M}}(K^{\pNum}_t)^2 \nu_t^\ell \wt \nu\upN_t  \\ & \quad  - 2 (1 - K\upN_t C) q_t^\ell K\upN_t (\nu_t^\ell -\frac{1}{\sqrt{M}} K\upN_t \wt \nu\upN_t )  \Big] \dd N_t.
\end{aligned}
\]%
Note that, we can write $Q\upN_t = \sum_{\ell =1}^\pNum (q_t^\ell)^2$. Furthermore, recalling the expression for empirical gain in \eqref{eq:defEmpGain} we can check that
\[
(K^{\pNum}_t)^2 C^2 Q\upN_t- 2 K\upN_t C Q\upN_t  + (K\upN_t)^2 V = - K\upN_t C Q\upN_t.
\]
The last two identities allow us to write
\begin{equation}\label{eq:covarMparticles}
\begin{aligned}
& \dd Q\upN_t =  \left (2 A Q\upN_t + G^2 - \frac{1}{M} G^2\right) \dd t 
 + 2 G \sum_{\ell = 1}^{\pNum} q_t^\ell \left(\dd \omega_t^\ell - \frac{1}{\sqrt{M}}\dd \wt \omega\upN_t \right) \\
& \qquad + \Big[ - K\upN_t C Q\upN_t - (K\upN_t)^2 V + \frac{1}{M}(K^{\pNum}_t)^2 \sum_{\ell=1}^{\pNum} (\nu_t^\ell)^2  \\
& \qquad  + \frac{1}{M}(K^{\pNum}_t)^2 (\wt \nu\upN_t)^2 - \frac{2}{M\sqrt{M}}(K^{\pNum}_t)^2 \sum_{\ell=1}^{\pNum}\nu_t^\ell \wt \nu\upN_t   \\
& \qquad - \frac{2}{M} (1 - K\upN_t C) K\upN_t \sum_{\ell=1}^{\pNum} q_t^\ell (\nu_t^\ell -\frac{1}{\sqrt{M}} K\upN_t \wt \nu\upN_t ) \Big] \dd N_t.
\end{aligned}
\end{equation}%
We next use this expression to compute the expectation of $\EE\expecof{Q\upN_t}$ with respect to the noise terms and, in addition, the sampling times of the Poisson counter.

\subsection{Convergence of expected covariance}
\begin{verArxiv}
Using Proposition~\ref{prop:generator}, we now turn our attention to computing the expectation of optimal covariance $P_t$ (with respect to sampling times), which was described by a piecewise deterministic Markov process in \eqref{eq:optVarPoisson}. Similarly, we will use Proposition~\ref{prop:generator} to describe the expectation of the empirical covariance (with respect to noise in the ensemble filter and the underlying sampling process), which is described by a jump diffusion process in \eqref{eq:covarMparticles}.
\end{verArxiv}
\begin{verMTNS}
We now use the generator based approach to compute the expectation of optimal covariance $P_t$ (with respect to sampling times), which was described by a piecewise deterministic Markov process in \eqref{eq:optVarPoisson}. We use similar technique to describe the expectation of the empirical covariance (with respect to noise in the ensemble filter and the underlying sampling process), which is described by a jump diffusion process in \eqref{eq:covarMparticles}.    
\end{verMTNS}
In what follows, we use the notation $\cP_t := \EE\expecof{P_t}$, $\cQ\upN_t = \EE\expecof{\cQ\upN_t}$, and similarly, for the corresponding injection gains, we use $\cK_t := \cP_t C(\cP_t C^2 + V)^{-1}$ and $\cK\upN_t := \cQ\upN_t C(\cQ\upN_t C^2 + V)^{-1}$.
\begin{verArxiv}
By applying Proposition~\ref{prop:generator} to the process \eqref{eq:optVarPoisson}, we first get    
\end{verArxiv}
\begin{verMTNS}
For the process \eqref{eq:optVarPoisson}, it can be shown that
\end{verMTNS}
\begin{equation}\label{eq:expectedCovOpt}
\begin{aligned}
\dot \cP_t &= 2 A \cP_t + G^2 - \lambda \cP_t^2 C^2 ( \cP_t C^2 + V)^{-1}  \\
&= (2 A -\lambda ) \cP_t + G^2 + \lambda \cP_t ( 1 - \cP_t C^2  (\cP_t C^2 + V)^{-1} ) \\
& = (2 A -\lambda ) \cP_t + G^2 + \lambda \cP_t \phi(\cP_t) 
\end{aligned}
\end{equation}
where the mapping $\cP_t \mapsto \phi(\cP_t)$ is defined as
\[
\phi(\cP_t) := (1 - \cP_t C^2 (\cP_t C^2 + V)^{-1} ).
\]
On the other hand, using the fact that $\nu_t^\ell$ and $\wt \nu\upN_t$ are mean zero processes, 
\begin{verMTNS}
it can be shown that for the process $Q\upN_t$ in \eqref{eq:covarMparticles}, we obtain
\end{verMTNS}
\begin{verArxiv}
the application of Proposition~\ref{prop:generator} to the process $Q\upN_t$ in \eqref{eq:covarMparticles} yields,
\end{verArxiv}
{\small 
\begin{equation}\label{eq:expectedCovEmp}
\begin{aligned}
\dot \cQ\upN_t & = 2 A \cQ\upN_t + G^2 - \frac{1}{M} G^2 - \lambda \cQ\upN_t C^2 (\cQ\upN_t C^2 + V)^{-1} \cQ\upN_t + \frac{\lambda}{M} (\cK\upN_t)^2 V \\
& = (2 A -\lambda) \cQ\upN_t + G^2 + \lambda \cQ\upN_t \phi(\cQ\upN_t) - \frac{1}{M} G^2 + \frac{\lambda}{M} (\cK\upN_t)^2  V .
\end{aligned}
\end{equation}
}%
To analyze the difference between optimal and empirical covariance, let us introduce $\cE_t\upN := \cQ_t\upN - \cP_t$, which satisfies the following equation:
\begin{equation}\label{eq:expectedErrorCov}
\begin{aligned}
\dot \cE_t\upN & = (2A - \lambda) \cE_t\upN + \lambda \cQ\upN_t \phi(\cQ\upN_t) - \lambda \cP_t \phi(\cP_t) - \frac{1}{M} G^2 + \frac{\lambda}{M} (\cK\upN_t)^2 V .
\end{aligned}
\end{equation}
The following lemma states some important properties that we need in the sequel. 
\begin{verArxiv}
    Its proof can be found in Appendix~\ref{app:boundSamplingRate}.
\end{verArxiv}
\begin{verMTNS}
    Its proof can be found in the extended version \citep{ref:TanYuf-24}.
\end{verMTNS}
\begin{lem}\label{lem:phiBounds}
For the equation~\eqref{eq:expectedCovOpt} and equation~\eqref{eq:expectedCovEmp}, we have the following relations:
\begin{enumerate}
\item There exists $0 < \gamma \le \ol\gamma < 1$ such that for each $t \ge 0$, we have
\[
\begin{aligned}
1 - \gamma \le \| \cP_tC^2 (\cP_tC^2 + V)^{-1} \| & < 1,\\
1 - \gamma \le \| \cQ\upN_t C^2 (\cQ\upN_t C^2 + V)^{-1} \| & < 1.
\end{aligned}
\]
Consequently, it holds that $\| \phi (\cP_t) \| \le  \gamma$ and $\| \phi(\cQ\upN_t) \| \le \gamma$.
\item The function $\phi$ has the property that for each $P_1$ and $P_2$, it satisfies
\[
\phi(P_1) P_1 - \phi(P_2)P_2 \le \phi(P_1)(P_1 - P_2) \phi(P_2).
\]
\end{enumerate}
\end{lem}

One can now use the statements in this lemma to get a bound on the evolution of $\cE_t\upN$. This is done by looking at the derivative of $(\cE\upN_t)^2$, that is,
\begin{align*}
\frac{d}{dt} (\cE_t\upN)^2 & \le 2 (2A - \lambda)(\cE_t\upN)^2 + 2 \lambda\gamma^2 (\cE_t\upN)^2 +\frac{2}{M} \|\cE_t\upN\| G^2 +\frac{2\lambda }{MC^2} V \|\cE_t\upN\|.
\end{align*}
In the last inequality, the quadratic term in $(\cE_t\upN)^2$ is negative if $\lambda - 2A > \lambda \gamma^2$,
or equivalently,
\begin{equation}\label{eq:lowerBndlambda}
\lambda > \frac{2A}{1- \gamma^2}.
\end{equation}
The ultimate bound on $\cE\upN_t$ now follows for every finite value of $\lambda$ satisfying \eqref{eq:lowerBndlambda}. In particular,
\[
\limsup_{t\to\infty} \| \cE\upN_t \| \le \frac{G^2 + \lambda V/C^2 } {M(\lambda(1-\gamma^2) - 2 A)}
\]
where we note that the denominator is strictly positive due to the bound on $\lambda$ in \eqref{eq:lowerBndlambda}.

\subsection{Convergence of expected mean}
The last step in the proof of Theorem~\ref{thm:main} is to analyze the difference between the expected values of the optimal mean and empirical mean. In this regard,  we let $\wh\cS\upN_t := \EE\expecof{\wh S\upN_t }$, and $\wh\cX_t = \EE\expecof{\wh X_t }$. 
\begin{verArxiv}
Applying Proposition~\ref{prop:generator}, we first observe that    
\end{verArxiv}
\begin{verMTNS}
We first observe that    
\end{verMTNS}
\[
\dot{\wh\cX_t} = A \wh\cX_t + \lambda \cK_t (\EE\expecof{y_t} - C\wh\cX_t ) 
\]
and
\[
\dot{\wh\cS\upN_t} = A \wh\cS\upN_t + \lambda \cK\upN_t (\EE\expecof{y_t} - C \wh\cS\upN_t)
\]
We now look at the asymptotic behavior of $\wh\chi\upN_t:=  \wh\cX_t - \wh \cS\upN_t$, which satisfies
\[
\begin{aligned}
\dot{\wh\chi}\upN_t  &= A \wh\chi\upN_t + \lambda \cK_t (y_t - C\wh\cX_t )- \lambda \cK\upN_t (\EE\expecof{y_t} - C \wh\cS\upN_t)  \\
& = (A - \lambda \cK\upN_tC) \wh\chi\upN_t + \lambda (\cK_t - \cK\upN_t )(\EE\expecof{y_t} - C \wh \cX_t).
\end{aligned}
\]
The bounds on asymptotic value of $\wh\chi\upN_t$ can now be established by analyzing the derivative of the following weighted quadratic function:
\[
\wh\chi\upN_t \mapsto (\cQ\upN_t)^{-1}(\wh\chi\upN_t)^2.
\]
We observe that
\[
\begin{aligned}
\frac{d}{dt} (\cQ\upN_t)^{-1} & = - (\cQ\upN_t)^{-1} \dot \cQ\upN_t (\cQ\upN_t)^{-1} \\
& = - 2 (A - \lambda \cK\upN_t C)(\cQ\upN_t)^{-1} - \lambda C^2(\cQ\upN_t C^2 + V )^{-1} \\ 
& \quad - (\cQ\upN_t)^{-2} \left( G^2 + \frac{1}{M} G^2 + \frac{\lambda}{M} \cK\upN_t V \cK\upN_t \right).
\end{aligned}
\]
Consequently, this last expression yields
\[
\begin{aligned}
\frac{d}{dt} \wh\chi\upN_t (\cQ\upN_t)^{-1} \wh\chi\upN_t & = - (\wh\chi\upN_t)^2  \Big[\lambda C^2(\cQ\upN_t C^2 + V )^{-1} \\ &  + (\cQ\upN_t)^{-2} \left( G^2 + \frac{1}{M} G^2 + \frac{\lambda}{M} \cK\upN_t V \cK\upN_t \right) \Big] \\ & + 2 \lambda (\cQ\upN_t)^{-1} \wh\chi\upN_t (\cK_t - \cK\upN_t )(\EE\expecof{y_t} - C \wh \cX_t).
\end{aligned}
\]
The right-hand side can be seen as a sum of two terms: the first one is negative definite and is quadratic in $\wh\chi\upN_t$; the second term is linear in $\wh\chi\upN_t$ but contains $(\cK_t - \cK\upN_t)$ which is converging to zero as $\pNum$ gets large. Using simple algebraic manipulations, and noting that $\cQ\upN_t$ is positive definite, it can be shown that $\wh\chi\upN_t$ converges to zero, as $t\to \infty$ and as $\pNum$ gets large. This completes the proof of Theorem~\ref{thm:main}.

\section{Conclusions}\label{sec:conclusions}
In this paper, we considered the analysis of vanilla ensemble filters for continuous-time linear Gaussian processes with Poisson-sampled discrete observations. The description of the particles is based on a continuous-discrete mean field process that contains flow equations with Wiener processes and jump equations with additive Gaussian noise in the modeled observation process. Using an extended generator, we compute the expectation of the empirical mean and empirical covariance of the particles with respect to the sampling process. The main result shows that for a sufficiently large value of the mean sampling rate, the empirical mean and covariance converge to the optimal mean and variance, respectively, as the number of particles tend to infinity.

The immediate extension of our work would be to generalize the proof of Theorem~\ref{thm:main} to systems of higher dimension. A better understanding of the bounds on the mean sampling rate needs is required. It is also interesting to investigate whether the tools developed in this paper can be used to study broader class of particle filters with stochastic noise.

{
\bibliography{biblio/Particles.bib,biblio/Filtering.bib,biblio/Books.bib,biblio/NewReferences.bib}
}

\begin{verArxiv}

\appendix

\section{Tools for Analysis}\label{sec:tools}
In this section, we discuss some analysis tools for stochastic differential equations that will be used in Section~\ref{sec:analysis} for the proof of Theorem~\ref{thm:main}. To do so, we consider the following differential equation for a stochastic process $\proc{x_t}_{t\ge 0}$ evolving in $\R^n$:
\begin{equation}\label{eq:sdeGeneral}
\dd x_t = f(x_t) \dd t + g(x_t) \dd \omega + h(x_t,\nu_{N_t}) \dd N_t
\end{equation}
where $N_t$ is a Poisson process with intensity $\lambda > 0$, $\omega$ is a standard Wiener process, and $\nu_{N_t}$ is a sequence of i.i.d. random variables with probability law $\mu$. 

To study the evolution of a function of the random process $\proc{x_t}_{t\ge 0}$, we make use of the Ito's chain rule. The reader may consult \citep[Chapter~II, Section 7]{protter2005stochastic} for detailed exposition on this topic. Here, the particular form we adopt is tailored for the differential equations appearing in earlier sections.

\begin{prop}[Ito's chain rule]\label{prop:chainRuleIto}
For a twice continuously differentiable function $\psi:\R^n \to \R$, it holds that
\begin{multline}\label{eq:itoChainRule}
\dd \psi (x_t)  = \Big[ \langle \nabla \psi(x_t), f(x_t) \rangle + \frac{1}{2}{\rm tr} \Big(\frac{\partial^2 \psi(x_t)}{\partial x} g(x_t) g(x_t)^\top \Big) \Big] \dd t \\ + \langle \nabla \psi(x_t), g(x_t) \rangle \dd \omega + \Big[ \psi(x_t +h(x,\nu_{N_t})) - \psi(x_t)\Big] \dd N_t.
\end{multline}
\end{prop}
Ito's chain rule describes the evolution of the function $\psi$ evaluated along the solution of the stochastic differential equation \eqref{eq:sdeGeneral}. However, to describe the evolution of the expectation of $\psi$ in differential form, we need to consider the extended generator as defined below:
\begin{defn}[Extended generator]
Given a real-valued function $\psi:\R^n \to \R$, the \emph{extended generator} of the process \( \proc[\big]{x_t}_{t\ge 0}\) described by \eqref{eq:sdeGeneral} is the linear operator \(\psi\mapsto \infgen\psi\) defined by
\begin{equation}
	\label{eq:defInfGen}
	\begin{aligned}
		 &\qquad \qquad \qquad \R^n \ni z \mapsto \infgen\psi(z) \in\R \\
		& \infgen\psi(z) \Let  \lim_{\varepsilon \downarrow 0} \frac{1}{\eps}\Bigl( \EE\cexpecof[\big]{ \psi \bigl( x(t+\eps)\bigr) \given x(t) = z} - \psi(z) \Bigr) .
	\end{aligned}
\end{equation}
\end{defn}
 
We obtain the expected value of $\psi$ by integrating the generator, which can be seen as a generalization of the classical Dynkin's formula:
\begin{equation}
	\label{e:Dynkin formula}
	\EE\expecof[\big]{\psi(x(t))} = \EE\expecof[\big]{\psi(x(0))} + \EE\expecof[\bigg]{\int_0^t \infgen \psi (x(s)) \, \dd s}. 
\end{equation}
For our purposes, it is useful to compute an explicit expression of the generator which can then be analyzed for studying the qualitative behavior of $\EE\expecof{\psi(x)}$. Several references in the literature derive generator equations for stochastic processes with jumps, see for example \citep{ref:HesTee-06} for a derivation in the context of Wiener process driven differential equations with renewal processes. 
\begin{prop}\label{prop:generator}
	If the sampling process \(\proc{N_t}_{t\ge 0}\) is Poisson with intensity \(\lambda > 0\), then the process \( \proc[\big]{x(t)}_{t\ge 0}\) described in \eqref{eq:sdeGeneral} is Markovian. Moreover, for any function $\R^n \ni z \mapsto \psi(z) \in \R $ with at most polynomial growth as \(\norm{z}\to+\infty\), we have
\begin{multline}\label{eq:formulaGen}
\cL \psi(z) = \langle \nabla \psi(z), f(z) \rangle + \frac{1}{2}{\rm tr} \Big(\frac{\partial^2 \psi(z)}{\partial z^2} g(z) g(z)^\top \Big) \\  + \lambda \Big[ \int \psi\bigl(z + h(z,\nu)\bigr) \mu(\dd\nu) - \psi(z) \Big].
\end{multline}
\end{prop}

\section{Proof of Lemma~\ref{lem:phiBounds}}
We provide a proof for the third item, which is independent of the first two items. To do so, we introduce the notation 
\[
K_1 := P_1 C (P_1C^2 + V)^{-1}; \quad K_2 := P_2 C (P_2C^2 + V)^{-1}.
\]
It can be easily checked that
\[
\phi(P_1) = 1 - K_1 C; \quad \phi(P_2) = 1 - K_2 C,
\]
and moreover, we have the identities
\[
\phi(P_1) P_1 C = K_1 V; \quad \phi(P_2) P_2 C = K_2 V.
\]
The last two equations result in
\begin{align*}
& \phi(P_1) P_1 - \phi(P_2)P_2  \\ 
& =  \phi(P_1) P_1 + K_1 V K_2 - K_1 V  K_2 - \phi(P_2)P_2  \\ 
& = \phi(P_1) P_1 + \phi(P_1) P_1 C K_2 - K_1 \phi(P_2) P_2 C  - \phi(P_2)P_2 \\
& = \phi(P_1) P_1 (1 - K_2 C) - (1 - K_1 C) \phi(P_2)P_2 \\
& = \phi(P_1) P_1 \phi(P_2) - \phi(P_1) \phi(P_2)P_2 \\
& = \phi(P_1)(P_1 - P_2) \phi(P_2)
\end{align*}
which is the desired result.

\section{Bounds on Sampling Rate}\label{app:boundSamplingRate}

It can be seen that the constant $\gamma$ used in the derivations actually depends on the mean sampling rate $\lambda$. Here we derive the bounds explicitly on $\lambda$. From the definition of $\gamma$, it follows that
\[
\gamma \le 1 - \frac{\cP_t C^2}{\cP_t C^2 + V} = \frac{V}{\cP_t C^2 + V}.
\]
On the other hand, it can be shown that
\[
\cP_t \ge \min \left\{P_0, \frac{G^2}{\lambda-2A} \right\}.
\]
and, for $G_M:=G\sqrt{(1-1/M)}$
\[
\cQ_t \ge \min \left\{Q_0^M, \frac{G_M^2}{\lambda-2A} \right\}.
\]
Letting $\wt \lambda := \lambda - 2A$, we get
\[
\gamma \le \max\left\{\frac{V}{\min\{P_0,Q_0^M\} C^2 + V},  \frac{V \wt \lambda}{G_M^2 C^2 + V \wt \lambda} \right\}\ .
\]
Consequently, the condition
\[
\wt \lambda > \lambda \gamma^2
\]
holds if 
\[
\wt \lambda > \max\left\{\frac{V^2}{(\min\{P_0,Q_0\} C^2 + V)^2},  \frac{V^2 \wt \lambda^2}{(G_M^2 C^2 + V \wt \lambda)^2} \right\} (\wt \lambda + 2A)\ .
\]
Using simple algebraic manipulations, the last inequality holds if and only if
\[
2 V (\lambda - 2A)(AV - G_M^2 C^2) < G_M^4 C^4,
\]
and
\[
\lambda > 2A \left(1-\frac{V^2}{(\min\{P_0,Q_0\} C^2 + V)^2}\right)^{-1}.
\]
The first inequality holds if $AV \le G_M^2C^2$, or the measurement noise covariance $V$ is relatively small.

\end{verArxiv}

\end{document}